%% file: BraunDevelin_EhrhartRootsNonNegativity.tex
\documentclass{amsart}
\usepackage{graphicx}
\usepackage{amsfonts, amsmath, amssymb, amsthm}
\usepackage{latexsym, color}

\newtheorem{theorem}{Theorem}[section]

\newtheorem{conjecture}[theorem]{Conjecture}
\newtheorem{question}[theorem]{Question}

\theoremstyle{definition}
\newtheorem{definition}[theorem]{Definition}
\newtheorem{example}[theorem]{Example}

\theoremstyle{remark}

\numberwithin{equation}{section}

\begin{document}

\title[Ehrhart Polynomial Roots]{Ehrhart Polynomial Roots and Stanley's Non-negativity Theorem}

\author{Benjamin Braun}

\address{Department of Mathematics, Washington University, St. Louis, MO}
\urladdr{http://math.wustl.edu/$\sim$bjbraun}
\email{bjbraun@math.wustl.edu}

\author{Mike Develin}

\address{American Institute of Mathematics, 360 Portage Ave., Palo Alto, 
CA 94306}
\urladdr{http://math.berkeley.edu/$\sim$develin/}
\email{develin@post.harvard.edu}


\date{October 4, 2006.}


\begin{abstract}
Stanley's non-negativity theorem is at the heart of many of the results in Ehrhart theory.  In this paper, we analyze the root behavior of general polynomials satisfying the conditions of Stanley's theorem and compare this to the known root behavior of Ehrhart polynomials.  We provide a possible counterexample to a conjecture of the second author, M. Beck, J. De Loera, J. Pfeifle, and R. Stanley, and contribute some experimental data as well.
\end{abstract}

\maketitle

Let $P$ be a convex polytope in $\mathbb{R}^n$ with vertices in $\mathbb{Z}^n$ and affine span of dimension $d$.  We will refer to such polytopes as \textit{lattice polytopes} and to elements of $\mathbb{Z}^n$ as \textit{lattice points}.  By a remarkable theorem due to E. Ehrhart, \cite{Ehrhart}, the number of lattice points in the $t^{th}$ dilate of $P$, for non-negative integers $t$, is given by a polynomial in $t$ of degree $d$ called the \textit{Ehrhart polynomial} of $P$.  In this paper we will investigate some differences between the root behavior of Ehrhart polynomials for elements of an arbitrary collection of polytopes of dimension less than or equal to $d$ and the root behavior of an arbitrary collection of polynomials of degree less than or equal to $d$ satisfying a certain non-negativity condition.

\section{Ehrhart Theory}

We will begin by reviewing basic facts about Ehrhart polynomials. It is well known, e.g. chapter 4 of \cite{StanleyVol1}, that for a polynomial $f$ of degree $d$ over the complex numbers there exist complex values $h_j^*$ so that \[\frac{\sum_{j=0}^{d}h_j^*x^j}{(1-x)^{d+1}}=\sum_{t\geq 0}f(t)x^t.\] \noindent Given a lattice polytope $P$, we denote its Ehrhart polynomial by $L_P(t)$ and let\[\mathrm{Ehr}_P(x)=\sum_{t\geq 0}L_P(t)x^t=\frac{\sum_{j=0}^dh_j^*x^j}{(1-x)^{d+1}}\] denote the \textit{Ehrhart series for $P$}.  There is a well-known relationship between this encoding of $L_P(t)$ and the polynomial itself, namely that $L_P(t)$ can be expressed as  \[L_P(t)=\sum_{j=0}^{d}h_j^*{t+d-j \choose d}.\]  This is easily seen by expanding the rational function as a formal power series.  

Thus, $\mathrm{Ehr}_P(x)$ encodes the change of coefficients for $L_P(t)$ corresponding to the change from the standard monomial basis to the basis \[B_d := \left\{ {t+d-j \choose d} : 0 \leq j \leq d \right\}.\]  \noindent It turns out that representing $g(t)$ in this way can be very profitable.  One of the most important known results about Ehrhart polynomials is the following theorem due to R. Stanley, known as Stanley's non-negativity theorem.

\begin{theorem} (see \cite{StanleyDecompositions} and \cite{BeckRobinsCCD}) If $P$ is a $d$-dimensional lattice polytope, then $(h_0^*,\ldots,h_d^*)\in (\mathbb{Z}_{\geq 0})^{d+1}$.
\end{theorem}

The non-negativity theorem is at the heart of much of what is presently known about roots and coefficients of Ehrhart polynomials.  However, not every polynomial with non-negative integer $h_j^*$'s is an Ehrhart polynomial, hence we make the following definition.

\begin{definition} A non-zero polynomial satisfying the condition that $$(h_0^*,\ldots,h_d^*) \in (\mathbb{R}_{\geq 0})^{d+1}$$ is called a \textit{Stanley non-negative, or SNN}, polynomial.
\end{definition}

In this paper we will see some of the differences between the root behavior of Ehrhart polynomials and of arbitrary SNN polynomials.

\section{Norm Bounds and Growth Rates}

In this section we review a norm bound on roots of SNN polynomials and some results and conjectures about growth rates of roots of SNN and Ehrhart polynomials.  In \cite{BDDPS}, it was shown that for polytopes of fixed dimension $d$ the roots of $L_P(t)$ are bounded above in norm by $1+(d+1)!$.  It was further suggested that this might be made polynomial in $d$.  In response, the first author proved the following.

\begin{theorem}\label{BraunRootThm} (see \cite{BraunRoots}) If $f$ is an SNN polynomial, then all the roots of $f$ lie inside the closed disc with center $\frac{-1}{2}$ and radius $d(d- \frac{1}{2})$.
\end{theorem}

The proof of this can be found in \cite{BraunRoots}, but we will find it useful to sketch the argument.  If we represent such an $f$ as \[f(t)=\sum_{j=0}^{d}h_j^*{t+d-j \choose d},\] then to evaluate $f$ at $t\in \mathbb{C}$ we take a non-negative linear combination of the $d+1$ points ${t+d-j \choose d}\in \mathbb{C}$, $0\leq j \leq d$.  If all of these points are in a common half space $H$ of $\mathbb{C}$ with zero on the boundary, then $f(t)\neq 0$.  Thus, one only needs to show that this is satisfied for $t\notin \{z:|z+\frac{1}{2}|\leq d(d-\frac{1}{2})\}$ to prove Theorem \ref{BraunRootThm}.

The above bound is essentially optimal, as the following theorem demonstrates.

\begin{theorem}\label{BeyHenkWillsThm}(see \cite{BeyHenkWills}) The polynomial $S_d(t)=\sum_{j=0}^{d}{t+d-j \choose d}$ is an Ehrhart polynomial whose roots all have real part $\frac{-1}{2}$.  Further, if $\alpha_d$ is the root of $S_d(t)$ of maximal norm, then \[\left| \alpha_d + \frac{1}{2}\right| = \frac{d(d+2)}{2\pi} + O(1)\] \noindent as $d\rightarrow \infty$.
\end{theorem}

Thus, a norm bound for Ehrhart polynomial roots cannot be better than quadratic in $d$.  It was suggested in \cite{BeyHenkWills} that $S_d(t)$ possesses the roots of maximal norm among all dimension $d$ polytopes with interior lattice points, and this was proved for $d=2,3$.  In response to the analogous question for SNN polynomials, we offer the following theorem and conjecture.

\begin{theorem}For the polynomial $M_d(t)={t+d \choose d} + {t \choose d}$, which is not an Ehrhart polynomial, if $\beta_d$ is the root of $M_d(t)$ of maximal norm, then \[|\beta_d + \frac{1}{2}| = \frac{d^2}{\pi} + O(1),\] as $d\rightarrow \infty$.
\end{theorem}

\begin{proof}
Note that the following closely follows the proof of Theorem \ref{BeyHenkWillsThm} given in \cite{BeyHenkWills}.

For any lattice polytope $P$, $h_d^*$ is equal to the number of interior lattice points in $P$.  If this is non-zero, then $h_1^*$, which records $L_P(t)-d+1$ when $P$ is $d$-dimensional, must also be non-zero.  As this condition is not satisfied by $M_d(t)$, $M_d(t)$ is not an Ehrhart polynomial.

By a result of Rodriguez-Villegas in \cite{Rodriguez-Villegas}, since the roots of the numerator of the generating function for $M_d(t)$ lie on the unit circle, all the roots  of $M_d(t)$ are on the line $x=\frac{-1}{2}$.  If $s=\frac{-1}{2}+bi$, $b\geq 0$, is a root of $M_d(t)$, then we have \begin{equation}\label{Root}(s+d)(s+d-1)\cdots(s+1)=-s(s-1)\cdots(s-d+1),\end{equation} as any root $s$ of $M_d(t)$ satisfies \[-{s+d \choose d}={s \choose d}.\]  Writing $s-j=s_j=r_je^{i\theta_j}$ and noting that $|s+j+1|=|s-j|$ implies $s+j+1=r_je^{i(\pi -\theta_j)}$, we can rewrite \eqref{Root} as \[(-1)^{d+1}=e^{i(2\theta_0 + \cdots +2\theta_{d-1})}.\]  We now substitute $\frac{\pi}{2} + \phi_j=\theta_j$, where $\phi_j \in (0,\frac{\pi}{2}]$.  This gives a new equation, \[-1=e^{i(2\phi_0 + \cdots +2\phi_{d-1})}.\]  Therefore, we must have, for some positive odd value of $k$, \[\frac{k\pi}{2} = \sum_0^{d-1} \phi_j.\]  By definition, $\cot \phi_j = \frac{b}{j+\frac{1}{2}}$ for $j=0,\ldots,d-1$.  Thus, $s$ is a root of $M_d(t)$ of maximal imaginary part if and only if
\begin{equation}\label{CotangentFormula}p_d(b):=\sum_{j=0}^{d-1}\cot^{-1}\left(\frac{b}{j+\frac{1}{2}}\right)=\frac{\pi}{2},\end{equation} as each $p_d(b)$ is a strictly decreasing function of $b$.  Say that $p_d(b_d)=\frac{\pi}{2}$.

For $x>1$, we have \[\cot^{-1}(x)=\tan^{-1}(\frac{1}{x})=\sum_{k=0}^\infty \frac{(-1)^k}{(2k+1)(x^{2k+1})}.\]  By truncating the Taylor series after the first and second term, we have $\frac{1}{x} > \cot^{-1}(x) > \frac{1}{x} - \frac{1}{3x^3}$.  Using this inequality on each summand in \eqref{CotangentFormula}, substituting $x=\frac{b}{j+\frac{1}{2}}$ for each $j$, and then using Faulhaber's formulas for the resulting $j^m$ terms, we have that for $b>d+\frac{1}{2}$, \[\frac{d^2}{2b} > p_d(b) > \frac{d^2}{2b} - \frac{d^4}{108b^3}.\]

Suppose now that $\hat{b} = \frac{d^2}{\pi} - \alpha$, where $\alpha$ is some large constant.  In that case we have \begin{equation}\label{Inequality}p_d(b) > \frac{d^2}{2\hat{b}} - \frac{d^4}{108\hat{b}^3} = \frac{\pi d^2}{2(d^2 - \pi\alpha)} - \frac{\pi^3 d^4}{108(d^2-\pi\alpha)^3}.\end{equation}  The limit of the right hand side of \eqref{Inequality} as $d$ increases is $\frac{\pi}{2}$, and for all large enough $d$ the right hand side is greater than $\frac{\pi}{2}$.  As each $p_d(b)$ is decreasing in $b$, we see that $b_d \geq \hat{b} = \frac{d^2}{\pi} - \alpha$.  As we also have $\frac{d^2}{\pi} > b_d$ for large $d$, we have our result.
\end{proof}

\begin{conjecture}\label{RootConjecture}The root of the polynomial $M_d(t)$ with largest norm, call it $\gamma_d$, has $|\mathrm{Im}(\gamma_d)|$ maximal among the imaginary parts of all roots of degree $d$ SNN polynomials.
\end{conjecture}

Experimental data for the roots of a large number of SNN polynomials of degree less than or equal to seven form the basis for this conjecture.  As a corollary of Theorem \ref{thm:3DBound} below, Conjecture \ref{RootConjecture} is true when $d=2$.


\section{General Bounds and the Vertical Strip Conjecture}

In this section we are interested in the following conjecture, due to the second author, M. Beck, J. De Loera, J. Pfeifle, and R. Stanley.

\begin{conjecture}(see \cite{BDDPS}) If $L_P(t)$ is an Ehrhart polynomial of degree $d$, then the roots $\alpha_i$ of $L_P(t)$ satisfy $-d\leq \mathrm{Re}(\alpha_i) \leq d-1$ for all $i$.
\end{conjecture}

We will refer to this as the \textit{Vertical Strip Conjecture}, for 
obvious reasons.  The original motivation for this claim was experimental 
data produced from polytopes of relatively low dimensions, along with the 
fact that all real roots of a degree $d$ SNN polynomial lie in the 
interval $[-d,d-1]$ (as shown in \cite{BDDPS}).  If this conjecture is 
true, one might hope that, like the norm bounds above, it actually holds 
for SNN polynomials.  One step in this direction is the following.

\begin{theorem}\label{ConeTheorem}
For $d\geq 2$, let $C_d^1$ (respectively $C_d^2$) be the open pointed cone in $\mathbb{C}$ with vertex $d-1$ (respectively $-d$), angular width $\frac{2\pi}{d}$, and bisecting ray $\left( d-1,\infty \right)$ (respectively $\left( -\infty,-d \right)$).  For any SNN polynomial $f$ of degree $d$, $C_d^1 \cup C_d^2$ does not contain a root of $f$.
\end{theorem}

\begin{proof} That the real roots of such polynomials are in the interval $[-d,d-1]$ was shown in \cite{BDDPS}.  We prove that the theorem holds for complex values in $C_d^1$ with positive imaginary part.  The proof is similar for the other cases.  Let $t\in C_d^1$.  Then the difference between the arguments of ${t+d-j \choose d}$ and ${t+d-j-1 \choose d}$ is less than $\frac{\pi}{d}$ for all $j$, as this is equal to the difference between the arguments of $(t+d-j)$ and $(t-j)$, both of which have argument less than $\frac{\pi}{d}$.  Therefore, the points ${t+d-j \choose d}$, $0\leq j \leq d$, lie in a common half plane with zero on the boundary.  A nonnegative, nonzero linear combination of such points cannot be zero.
\end{proof}

Note that the Vertical Strip Conjecture follows from Theorem \ref{ConeTheorem} for $d=2$, as in this 
case the angular widths of our cones are $\pi$. We can extend this as follows.

\begin{theorem}\label{thm:VSC34}
For any SNN polynomial of degree $3$ or $4$, the Vertical Strip Conjecture holds.
\end{theorem}

\begin{proof}

If $d = 3$ or $4$, we show that for every complex number $z$ lying outside the vertical strip $\{z:-d\leq \mathrm{Re}(z) \leq d-1\}$, the numbers ${z+d-j \choose d}$, $0\le j\le d$, all lie in a half-plane with zero on the boundary.  This is tantamount to showing that the angles $A(j)$ formed between the vectors $z-j$ and $z+d-j$, which are the angle differences between successive ${z+d-j\choose d}$, sum to less than $\pi$ for $0\le j\le d-1$.  The situations for the real part of $z$ being less than $d$ and greater than $d-1$ are symmetric (interchanging $j$ with $d-1-j$), so assume that the real part of $z$ is greater than $d-1$.  As these are real polynomials, complex roots will occur in conjugate pairs, hence we may assume that all relevant vectors lie in the first quadrant and therefore each $A(j)$ is increased by decreasing the real part of $z$.  Hence, we need to show that the angle sum $\sum_jA(j)$ is at most $\pi$ when the real part of $z$ is $d-1$.

Let $z = (d-1) + ki$. An application of the Law of Cosines to the triangle with vertices 0, $z-j$, and 
$z+d-j$ yields the following formula, where $r = (d-1)-j$ (the $x$-coordinate of $z-j$) and $s = (2d-1-j)$ 
(the $x$-coordinate of $z+d-j$):
\[
\text{cos}^2 A(j) = 1 - \frac{d^2 k^2}{(k^4 + (r^2+s^2) k^2+ r^2s^2)}
\]

Since $d$ is a constant, applying the AM/GM inequality to the fraction yields that this quantity is 
minimized when $k^2 = rs$, and increases (thus A(j) decreases)  monotonically in both directions as $k$ 
gets further from this point.

Now, we move on to the specific applications for $d=3$ and $d=4$. For $d=3$, simply applying the above 
formula yields (by numerical computation):

\[
(A(0), A(1), A(2)) \le (0.45, 0.65, 1.58);
\]
the sum of these is less than $\pi$, so the relevant ${z+d-j\choose d}$ all lie in a half-plane.

For $d=4$, we need to consider cases. Using the monotonicity results 
obtained above, we divide into the 
cases $k^2 \le 2$ and $k^2\ge 2$. For $k^2\le 2$ we have:
\[
(A(0), A(1), A(2), A(3))\le (0.24, 0.38, 0.70, 1.58),
\]
while for $k^2 \ge 2$ we have:
\[
(A(0), A(1), A(2), A(3))\le (0.41, 0.53, 0.73, 1.23);
\]
in both cases the sum of these angles is less than $\pi$, and so in both cases the relevant 
${z+d-j\choose d}$ again all lie in a half-plane. This completes the proof that all zeroes lie in the 
vertical strip for $d=3, 4.$
\end{proof}

As an aside, the above results can be extended further to show that for $d=3$ or $4$, when a potential zero is not purely real, its real part must in fact satisfy tighter bounds than simply being between $-d$ and $d-1$. Indeed, the locus of all possible zeroes off the real axis is bounded by curves given by the equations $\sum A(j) = \pi$.  We have shown that these curves lie entirely inside the vertical strip for $d=3, 4$.  These curves will be discussed more in the final section.

While the above results are promising, in general the vertical strip conjecture is not satisfied by SNN 
polynomials. To produce examples illustrating this, we pick a desired root $z$ for which the numbers ${z+d-j \choose d}$ do not lie in a complex half-plane for a desired degree $d$. We can then produce a positive linear combination of these numbers which is equal to zero, and the polynomial encoded by these coefficients will have a root equal to $z$. For $d$ large enough, we can find such $z$ which lie outside the vertical strip, as the following examples demonstrate.

\begin{example}\label{SNNex}
The polynomial \[f(t) = {t+5\choose 5} + 33 {t\choose 5}\] satisfies $f(z) = 0$ for $z\approx 4.00019 + 3.00963i$.
\end{example}

We obtained this polynomial by noting that for $z = 4+3i$, the sum of the $A(j)$ exceeds $\pi$. Therefore, there is 
a SNN polynomial with root $4+3i$, and by changing the coefficients 
slightly, we obtain a polynomial with real part 
strictly greater than 4.

However, this polynomial is not an Ehrhart polynomial. In particular, it 
does not satisfy the list of inequalities satisfied by the coefficients of 
Ehrhart polynomials given in~\cite{BDDPS}. Furthermore:

\begin{theorem}
No SNN polynomial of degree 5 with a root outside the vertical strip can 
be an Ehrhart polynomial. Therefore, the Vertical Strip Conjecture holds 
for Ehrhart polynomials (though not SNN polynomials) for $d=5$.
\end{theorem}

\begin{proof}
We give a sketch of the proof. Using methods similar to those in the proof 
of Theorem~\ref{thm:VSC34}, one can show that the angle sum can only 
barely exceed $\pi$ (it is bounded by 3.17.) In particular, for any $d$ 
outside the vertical strip but inside the SNN locus, the vectors $v_j := {z+d-j\choose d}$ all lie in 
the same half-plane for $j\in \{0, \ldots, 4\}$ and also for $j\in \{1, \ldots, 5\}$.

This means that any positive dependence among these vectors is a sum of 
positive dependences among $\{v_0, v_j, v_5\}$ for $j\in [4]$. However, 
for all relevant $z$ (outside the vertical strip but inside the SNN 
locus), the magnitude of $v_5$ turns out to be much smaller than 
the magnitude of $v_0$ and $v_1$; in particular, in each of these dependences, 
as in Example~\ref{SNNex}, the coefficient of $v_5$ must be much larger 
than the coefficients of $v_0$ and $v_1$ (if $j=1$.)

This implies that $h_5^* > h_0^* + h_1^*$, which violates the inequality 
$h_5^*\le 
h_0^* + h_1^*$ given in~\cite{BDDPS}.
\end{proof}

However, by bumping up the degree, we can find candidate 
Ehrhart polynomials that may or may not be actual Ehrhart polynomials of 
polytopes.

\begin{example}\label{PossibleEhrhart}

Consider the polynomial $g(t)$ of degree $26$ with 
\[\left(h_0^*,\ldots,h_{26}^*\right)=\left(1,2,3,4,6,10,16,27,43,69,112,181,293,473,762,0,\ldots,0\right). \] 
Numerical approximation produces \[26.47331467-28.51231239i\] as a root of $g(t)$.

\end{example}

For large $d$, there are a large number of SNN polynomials with roots outside the vertical strip; any point for 
which the angle sum $\sum A(j)$ is larger than $\pi$ is the root of some SNN polynomial. If we pick a point close 
to the vertical strip, we will have a large number of SNN polynomials which have it as a root (an SNN polynomial is 
produced by any positive linear dependence among $d$ vectors which for large $d$ have roughly evenly spaced 
arguments). It seems as though one of these must be an Ehrhart polynomial, though of course discerning the 
Ehrhart-ness of degree 26 polynomials (which come from dimension 26 polytopes) with relatively large $h^*$-vectors is 
certainly a nontrivial task. Indeed, the following question is concrete and unresolved.

\begin{question}
Is the polynomial $g(t)$ from Example \ref{PossibleEhrhart} an Ehrhart polynomial?
\end{question}

The answer would be interesting either way; if it is an Ehrhart polynomial, it is a counterexample to the Vertical 
Strip Conjecture, and if not, new methods will need to be developed to verify this.

\section{Experimental Results and Bounds in Low Dimension}

For Ehrhart polynomials of degree $2$, there are very tight known restrictions on the location of the roots.

\begin{theorem}{\rm (Beck, et al, see \cite{BDDPS}.)} The roots of the Ehrhart polynomial of any lattice $2$-polytope are contained in \[\left\{ -2,-1,-\frac{2}{3} \right\} \cup\left\{x+iy\in \mathbb{C}:-\frac{1}{2}\leq x < 0, |y|\leq \frac{\sqrt{15}}{6}\right\}.\]

\end{theorem}

For degree $2$ SNN polynomials, we can produce very similar root restrictions.  As mentioned after the proof of Theorem \ref{thm:VSC34}, the zeros of SNN polynomials in fixed degree are contained in regions bounded by curves given by $\sum{A(j)}=\pi$.  By analyzing this curve in detail, we obtain the following.

\begin{theorem}\label{thm:3DBound}
The roots of any degree $2$ SNN polynomial are contained in \[[-3,2]\cup 
\left\{x+iy \in \mathbb{C}: y^2 \leq -x^2 -x +\frac{1}{2}\right\}.\]
\end{theorem}

\begin{proof}
Let $z = x+iy$. Suppose that $z$ is not pure real; consider the diagram in 
Figure~\ref{cosinelaw}, which 
shows the relevant points $z-1, z, z+1, z+2$. Then we need to determine 
the locus of all $z$ such that the angles $\theta$ and $\phi$ are 
complementary. The situation is obviously symmetric about the real axis 
and about $x = -1/2$, so we assume that $x>-1/2, y > 0$. 
In this picture, we must have $\theta > \phi$: the derivative of 
$\text{tan}^{-1}\,x$ is a 
decreasing function of $|x|$, the uncommon parts of $\theta$ and $\phi$ 
are the integral of this derivative over an interval of size 1, and the 
one which is part of $\phi$ is strictly further from the $y$-axis (since 
the real part of $z$ is greater than -1/2.) It then follows that $\theta + 
\phi > \pi$ if and only if $\text{sin } \theta < \text{sin } \phi$, or, 
since both sines are positive, if and only if $\text{sin}^2\,\theta < 
\text{sin}^2\,\phi$.

\begin{figure}[ht]
\begin{center}\input{thetaphi.pstex_t}\end{center}
\caption{\label{cosinelaw}}
\end{figure}

We use the Law of Sines, which states that the area of a triangle is 
one-half times the product of two adjacent sides times the sine of the 
included angle, on the two triangles with included angles $\theta$ and 
$\phi$. These triangles both have area equal to $y$, as their height is 
$y$ and their base is $2$. For the triangle with angle $\theta$, we obtain:

\begin{eqnarray*}
b & = & \frac{1}{2}\sqrt{(x-1)^2 + y^2}\sqrt{(x+1)^2 + y^2}\:
\text{sin }\theta, \text{ or} \\
\text{sin}^2\,\theta & = & \frac{4y^2}{((x-1)^2 + y^2))((x+1)^2 + y^2)}.\\
\end{eqnarray*}

Similarly, we obtain

\[
\text{sin}^2\,\phi = 
\frac{4y^2}{(x^2+y^2)((x+2)^2 + y^2)}.
\]

Comparing these formulas, we algebraically manipulate:

\begin{eqnarray*}
\frac{4y^2}{((x-1)^2 + y^2))((x+1)^2 + y^2)} &<& 
\frac{4y^2}{(x^2+y^2)((x+2)^2 + y^2)}
 \\
((x-1)^2 + y^2)((x+1)^2 + y^2)) &>& (x^2+y^2)((x+2)^2+y^2) \\
x^4 - 2x^2 + 1 + y^2 (2x^2 + 2) + y^4 &>& x^4 + 4(x^3 + x^2) + y^2 (2x^2 + 
4x+4) 
+ y^4 \\
y^2 (4x+2) + 4x^3 + 6x^2 - 1 &<& 0 \\
y^2 &<& \frac{1 - 4x^3 - 6x^2}{4x+2} \\
\end{eqnarray*}
where the last step is because $4x+2$ is positive. This reduces to $y^2 < 
-x^2 - x + 1/2$, which is the equation of the indicated circle.

Finally, if $z$ is pure real, the result follows immediately from 
Theorem~\ref{ConeTheorem}.
\end{proof}

There are also very tight restrictions on the location of roots of Ehrhart polynomials of degree $3$, as the following theorem shows.

\begin{theorem}{\rm (Bey, et al, see \cite{BeyHenkWills}.)} The roots of the Ehrhart polynomial of any lattice $3$-polytope are contained in \[[-3,1]\cup\left\{x+iy \in \mathbb{C}: -1\leq x < 1, x^2 + y^2 \leq 3\right\}.\]
\end{theorem}

\begin{figure}[!ht]
\centering
\includegraphics{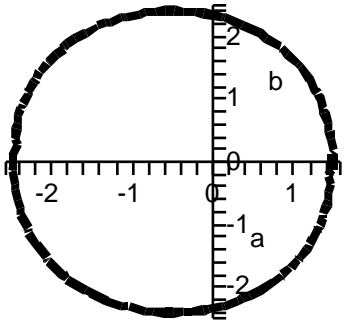}
\caption{The curve given by $\sum A(j) = \pi$ for degree $3$ SNN polynomials.}
\label{fig:RootBound}
\end{figure}

The curve $\sum{A(j)}=\pi$ for degree $3$ SNN polynomials is shown in Figure \ref{fig:RootBound}.  In degrees $2$ and $3$, note that the known restrictions for Ehrhart polynomial roots and SNN polynomial roots differ primarily in the restrictions on their real parts.

\begin{figure}[!ht]
\centering
\includegraphics{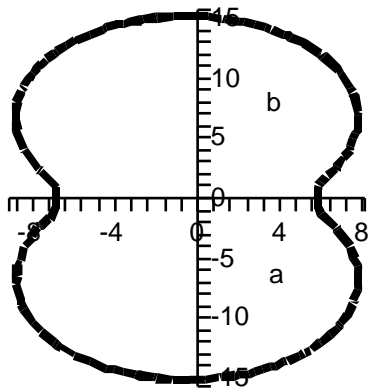}
\caption{The curve given by $\sum A(j) = \pi$ for degree $7$ SNN polynomials.}
\label{fig:7DRootBound}
\end{figure}

\begin{figure}[!ht]
\centering
\includegraphics{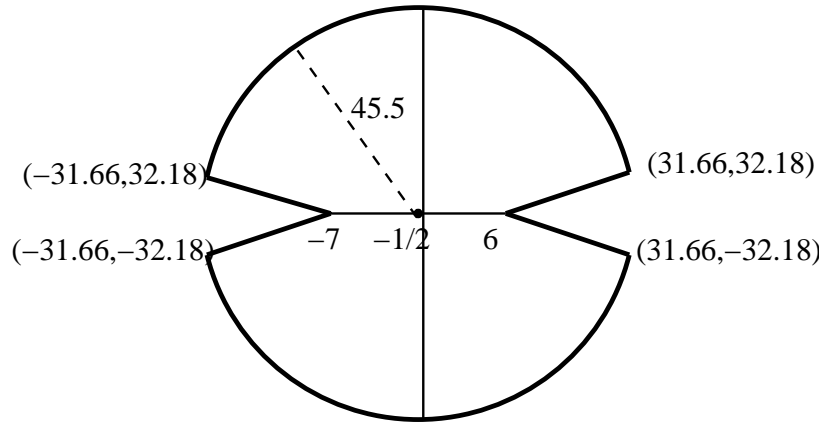}
\caption{The region containing the roots of degree $7$ SNN polynomials from Theorems \ref{BraunRootThm} and \ref{ConeTheorem}.}
\label{fig:7DThmBound}
\end{figure}

\begin{figure}[!ht]
\centering
\includegraphics{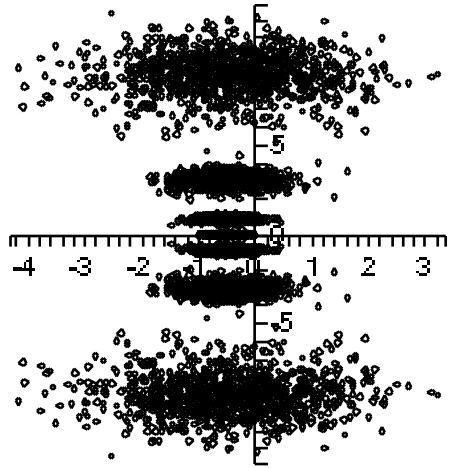}
\caption{The roots of $1000$ random degree $7$ SNN polynomials.}
\label{fig:7DRootPlot}
\end{figure}

Unfortunately, for higher degrees the bounding curves given by $\sum{A(j)}=\pi$ are more complicated.  Figure \ref{fig:7DRootBound} shows this curve for degree $7$ SNN polynomials.  It is interesting to compare this to Figure \ref{fig:7DRootPlot}, a plot of the roots of $1000$ random degree $7$ SNN polynomials, and to Figure \ref{fig:7DThmBound}, an approximation of the region containing the roots of degree $7$ SNN polynomials as determined by Theorems \ref{BraunRootThm} and \ref{ConeTheorem}.  All three of these pictures have roughly the same ``shape,'' though the random root plot is contained in a much smaller region than that bounded by the curve $\sum{A(j)}=\pi$, which is a smaller region than that given by Theorems \ref{BraunRootThm} and \ref{ConeTheorem}.

As we have seen, both Ehrhart and SNN polynomials have roots growing quadratically in norm as the degree grows, with the maximal roots of the extremal candidates for each class differing in norm by roughly a factor of two.  It would be interesting to see if there continues to be significant differences in the restrictions on the real parts of these roots in higher dimensions and if those differences can be quantitatively analyzed in a similar fashion.

\section{Acknowledgements}
Thanks to Martin Henk for his help regarding the proof of Theorem 
\ref{BeyHenkWillsThm} in \cite{BeyHenkWills}.  All graphics and numerical 
approximations in this paper were produced using the MAPLE computer 
algebra system. Part of this work was conducted while at the 2006 
AMS-IMS-SIAM Joint Summer Research Conference on Integer Points in 
Polyhedra in Snowbird, UT; we would like to thank the organizers for their 
support.

\break

\bibliographystyle{plain}
\bibliography{Braun}

\end{document}

%% file: thetaphi.pstex_t
\begin{picture}(0,0)%
\includegraphics{thetaphi.pstex}%
\end{picture}%
\setlength{\unitlength}{3947sp}%
\begingroup\makeatletter\ifx\SetFigFont\undefined%
\gdef\SetFigFont#1#2#3#4#5{%
  \reset@font\fontsize{#1}{#2pt}%
  \fontfamily{#3}\fontseries{#4}\fontshape{#5}%
  \selectfont}%
\fi\endgroup%
\begin{picture}(4076,1504)(4051,-2783)
\put(5851,-1711){\makebox(0,0)[lb]{\smash{{\SetFigFont{12}{14.4}{\familydefault}{\mddefault}{\updefault}{\color[rgb]{0,0,0}$\phi$}%
}}}}
\put(4801,-2011){\makebox(0,0)[lb]{\smash{{\SetFigFont{12}{14.4}{\familydefault}{\mddefault}{\updefault}{\color[rgb]{0,0,0}$\theta$}%
}}}}
\put(4051,-1411){\makebox(0,0)[lb]{\smash{{\SetFigFont{12}{14.4}{\rmdefault}{\mddefault}{\updefault}{\color[rgb]{0,0,0}$z-1$}%
}}}}
\put(5251,-1411){\makebox(0,0)[lb]{\smash{{\SetFigFont{12}{14.4}{\rmdefault}{\mddefault}{\updefault}{\color[rgb]{0,0,0}$z$}%
}}}}
\put(6451,-1411){\makebox(0,0)[lb]{\smash{{\SetFigFont{12}{14.4}{\rmdefault}{\mddefault}{\updefault}{\color[rgb]{0,0,0}$z+1$}%
}}}}
\put(7651,-1411){\makebox(0,0)[lb]{\smash{{\SetFigFont{12}{14.4}{\rmdefault}{\mddefault}{\updefault}{\color[rgb]{0,0,0}$z+2$}%
}}}}
\end{picture}%